\documentclass[brochure,12pt]{bourbaki}
\usepackage[matrix,arrow]{xy}
\usepackage{amssymb,amsfonts,amsmath,footnote}
\usepackage[francais]{babel}

\date{Mars 2013}
\bbkannee{65\`eme ann\'ee, 2012-2013}
\bbknumero{1067}
\title{Le mouvement brownien branchant vu depuis une extr\'emit\'e}
\subtitle{d'apr\`es Arguin-Bovier-Kistler et A\"id\'ekon-Berestycki-Brunet-Shi
}
\author{Jean-Baptiste Gou\'er\'e}
\address{Universit\'e d'Orl\'eans\\
UFR Sciences\\
MAPMO (UMR 6628 du CNRS)\\
F\'ed\'eration Denis Poisson\\
B.P. 6759\\
F--45067 Orl\'eans Cedex 2}
\email{Jean-Baptiste.Gouere@univ-orleans.fr}

\usepackage[utf8x]{inputenc}

\newcommand{\R}{\mathbb{R}}

\newcommand{\cF}{\mathcal{F}}
\newcommand{\cN}{\mathcal{N}}
\newcommand{\cNlineaire}{\mathcal{N}_{\ell}}
\newcommand{\cP}{\mathcal{P}}
\newcommand{\cDABK}{\mathcal{D}}
\newcommand{\cQ}{\mathcal{Q}}
\newcommand{\cC}{\mathcal{C}}
\newcommand{\cH}{\mathcal{H}}

\def\ABK{Arguin-Bovier-Kistler }
\def\ABBS{A\"id\'ekon-Berestycki-Brunet-Shi }
\def\ABetK{Arguin, Bovier et Kistler }
\def\ABBetS{A\"id\'ekon, Berestycki, Brunet et Shi }
\def\med{\hbox{med}}
\def\Cmed{C_{\tiny \med}}

\begin{document}
\maketitle


\section{Introduction}

Une particule est plac\'ee en l'origine de la droite r\'eelle \`a l'instant initial.
Cette particule se d\'eplace en suivant un mouvement brownien standard.
Au bout d'un temps ind\'ependant de loi exponentielle de moyenne $1$, cette particule se divise alors en $2$ particules.
Les deux particules ainsi obtenues \'evoluent ind\'ependamment et de la m\^eme mani\`ere que la premi\` ere particule : d\'eplacement brownien puis division.

Notons $N(t)$ le nombre de particules pr\'esentes \`a l'instant $t \ge 0$.
Remarquons que $N(t)$ est une variable aléatoire d'espérance $\exp(t)$.
Notons $X_1(t) \ge \cdots \ge X_{N(t)}(t)$ les positions de ces particules ordonnées par ordre décroissant.
Le comportement asymptotique du maximum $X_1(t)$ du mouvement brownien branchant a été l'objet de nombreux travaux de natures probabilistes et analytiques.
Cette popularité s'explique en partie par l'existence de liens très étroits entre ce modèle et la famille des équations aux dérivées partielles de
Fisher ou Kolmogorov-Petrovsky-Piscounov (F-KPP).
Ce lien a été mis en lumière par McKean \cite{McKean}.
Explicitons-le sur un exemple.
Notons $v(t,\cdot)$ la fonction de répartition du maximum du nuage de particules à l'instant $t \ge 0$ :
\begin{equation} \label{e:vintro}
v(t,x)=P(X_1(t) \le x).
\end{equation}
On vérifie facilement que cette fonction est solution de l'équation F-KPP suivante :
\begin{equation} \label{e:FKPPintro}
\frac{\partial v}{\partial t}
=
\frac12 \frac{\partial^2 v}{\partial^2 x} + v^2-v.
\end{equation}

Plus récemment, des résultats sur le comportement asymptotique de l'ensemble du nuage de particules vu depuis son maximum ont été obtenus
indépendamment par \ABK  \cite{ABK1,ABK2,ABK3,ABK4,ABK5} et \ABBS \cite{ABBS}.
Cet exposé est consacré à la présentation de ces résultats.
Nous nous concentrerons essentiellement sur les articles \cite{ABK3} et \cite{ABBS} et, dans une moindre mesure, \cite{ABK1}.
Ces travaux reposent notamment sur ceux de Bramson \cite{Bramson1,Bramson2}, Chauvin-Rouault \cite{Chauvin-Rouault-88,Chauvin-Rouault-90},
Lalley-Sellke \cite{Lalley-Sellke} et McKean \cite{McKean}.
Ils s'inspirent également en partie des travaux de Brunet-Derrida \cite{Brunet-Derrida-EPL,Brunet-Derrida-JSP}.

Plusieurs des propriétés du mouvement brownien branchant vu depuis un extremum sont conjecturées dans d'autres modèles,
notamment le champ libre gaussien en dimension $2$ 
\cite{Bolthausen-Deuschel-Giacomin, Bolthausen-Deuschel-Zeitouni, Bramson-Ding-Zeitouni,Bramson-Zeitouni,Ding-Zeitouni},
le temps de recouvrement de graphes par des marches aléatoires \cite{Dembo,Dembo-et-al}
et, plus généralement, les champs gaussiens exhibant des corrélations logarithmiques  \cite{Arguin-Zindy,Carpentier-et-al,Duplantier-al,Fyodorov-Bouchaud}.
Nous ne développerons pas ces aspects dans cet exposé.

\bigskip

Je remercie \'Elie Aïdékon, Louis-Pierre Arguin, Julien Berestycki, Anton Bovier, \'Eric Brunet, Nicola Kistler et Zhan Shi pour leurs réponses à mes questions
et pour leurs commentaires sur le manuscrit.

\section{Comportement asymptotique du maximum}

\subsection{Résultats}

Nous donnons principalement dans cette partie des résultats classiques sur le comportement asymptotique du maximum $X_1(t)$.

Notons $\med(t)$ la médiane du maximum $X_1(t)$ à l'instant $t$.
Notons $v(t,\cdot)$ sa fonction de répartition.
Elle est définie par \eqref{e:vintro} et elle vérifie l'équation F-KPP \eqref{e:FKPPintro}.
Les résultats de Kolmogorov-Petrovsky-Piscounov \cite{KPP} affirment l'existence d'une variable aléatoire non triviale $W$ et d'un terme de centralisation
$m(t) \sim \sqrt{2}$ tels que :
\begin{equation}\label{e:W}
X_1(t)-m(t) \to W \hbox{ en loi.}
\end{equation}
La fonction de répartition $w$ de $W$ est solution de 
$$
\frac12 w'' + \sqrt{2}w' + w^2-w = 0.
$$
Les solutions non triviales de cette équation sont uniques à translation près.
Il existe une constante $C_w>0$ (Bramson \cite{Bramson2}) telle que : 
\begin{equation} \label{e:asymptotiquew}
1-w(x) = P(W > x) \sim_{x \to +\infty} C_wxe^{-\sqrt{2}x}.
\end{equation}

Bramson \cite{Bramson2} a établi que l'on pouvait prendre pour $m$ la fonction définie par :
\begin{equation}\label{e:mtintro}
m(t) = \sqrt{2}t  - \frac{3}{2\sqrt{2}} \ln(t).
\end{equation}
Ainsi, il existence une constante $\Cmed$ telle que :
\begin{equation}\label{e:bramson-mediane}
\med(t) = \sqrt{2}t  - \frac{3}{2\sqrt{2}} \ln(t) + \Cmed + o(1).
\end{equation}

On définit une martingale (pour la filtration naturelle du mouvement brownien branchant) en posant, pour tout $t \ge 0$ :
\begin{equation}\label{e:Z}
Z(t) = \sum_{k=1}^{N(t)} \big(\sqrt{2}t-X_k(t)\big) e^{-\sqrt{2}\big(\sqrt{2}t-X_k(t)\big)}.
\end{equation}
Lalley et Sellke \cite{Lalley-Sellke} ont établi la convergence presque sûre de $Z(t)$ vers une variable aléatoire finie et strictement positive $Z$ ainsi
que la convergence presque sûre suivante :
\begin{equation} \label{e:LS}
\lim_{s \to \infty} \lim_{t \to \infty} P(X_1(t+s)-m(t+s) \le x | \cF_s) = \exp(-C_wZ e^{-\sqrt{2}x}) \hbox{ p.s.}
\end{equation}
On en déduit la représentation intégrale suivante :
\begin{equation}\label{e:representationw}
w(x) = P(W \le x) = E\left(\exp\left(-C_wZ e^{-\sqrt{2}x}\right)\right).
\end{equation}
\'Ecrivons
$$
\exp(-C_wZ e^{-\sqrt{2}x}) =  \exp(-e^{-\sqrt{2}(x - 2^{-1/2}\ln(C_wZ))}). 
$$
Conditionnellement à $Z$, c'est la fonction de répartition d'une distribution de Gumbel.
La convergence \eqref{e:LS} peut par conséquent s'interpréter de la manière suivante.
La variable aléatoire $X_1(t)-m(t)$ converge vers la somme de deux termes : 
le terme $2^{-1/2}\ln(C_wZ)$ qui provient de l'histoire du processus avant sa stabilisation en loi ;
un terme de fluctuation aléatoire qui suit une loi de Gumbel.
Plus précisément, Lalley et Sellke ont conjecturé la convergence en loi et la convergence en moyenne ergodique 
de l'ensemble du processus vu depuis $m(t)+2^{-1/2}\ln(C_wZ)$.
\ABK et \ABBS ont établi indépendemment ces conjectures et ont donné deux descriptions différentes du processus limite.

Signalons pour conclure cette partie qu'un résultat récent de Roberts \cite{Roberts} décrit le comportement presque sûr du maximum.
On a :
$$
\liminf_{t \to \infty} \frac{X_1(t) - \sqrt{2}t}{\ln(t)} \to  -\frac{3}{2\sqrt{2}} \hbox{ presque sûrement}
$$
mais
$$
\limsup_{t \to \infty} \frac{X_1(t) - \sqrt{2}t}{\ln(t)} \to  -\frac{1}{2\sqrt{2}} \hbox{ presque sûrement.}
$$
Ces résultats sont les analogues en temps continus de résultats obtenus peu auparavant par Hu et Shi \cite{Hu-Shi}.

\subsection{Quelques arguments}

Nous donnons dans cette partie les principaux arguments menant au résultat suivant :
\begin{equation}\label{e:bramson-facile}
\med(t) = \sqrt{2}t  - \frac{3}{2\sqrt{2}} \ln(t) + O(1).
\end{equation}
C'est une version faible de \eqref{e:bramson-mediane}.
Elle a été établie pour la première fois par Bramson \cite{Bramson1}.
Une preuve courte a été fournie par Roberts \cite{Roberts}.
Cette partie s'inspire en particulier de ce dernier article et de l'article de Addario-Berry et Reed \cite{Addario-Berry-Reed}.

\bigskip

Commençons par étudier le cas élémentaire où les particules sont indépendantes.
Plus précisément, donnons-nous, conditionnellement à $N(t)$, $N(t)$ particules de positions gaussiennes indépendantes centrées et de variance $t$.
Notons $X_1^*(t) \ge X_2^*(t) \ge \cdots \ge X_{N(t)}^*(t)$ les positions ordonnées des particules.
Désignons par $\med^*(t)$ la médiane de $X_1^*(t)$.
La médiane admet le développement asymptotique suivant :
$$
\med^*(t) = \sqrt{2}t  - \frac{1}{2\sqrt{2}} \ln(t) + O(1).
$$
Ce résultat est élémentaire.
Il peut par exemple être obtenu par des considérations sur les deux premiers moments de $N^*_{q(t)}(t)$,
le nombre de particules au-dessus de $q(t)$ à l'instant $t$ :
$$
N^*_{q(t)}(t) = \sum_{k=1}^{N(t)} 1_{X_k^*(t) \ge q(t)}.
$$
On a :
$$
\frac{E(N^*_{q(t)}(t))^2}{E(N^*_{q(t)}(t)^2)} \le P(N^*_{q(t)}(t) \ge 1) \le E(N^*_{q(t)}(t)).
$$
L'indépendance entre les particules permet de comparer utilement le second moment et le carré du premier moment.
La détermination, à une constante additive près, de la médiane $\med^*(t)$ est alors essentiellement ramenée à 
la détermination du réel $q(t)$ pour lequel le premier moment est d'ordre $1$.
Prenons $q(t)$ de la forme $\sqrt{2}t-a(t)$ où $a(t)$ est négligeable devant $\sqrt{t}$.
On a :
\begin{equation} \label{e:m1indep}
E \big(N^*_{q(t)}(t) \big) \sim \frac{1}{2\sqrt{\pi t}} e^{t-\frac{q(t)^2}{2t}} \sim \frac{1}{2\sqrt{\pi t}} e^{a(t)\sqrt{2}}.
\end{equation}
Cette quantité est d'ordre $1$ pour tout $t$ grand lorsque 
$$
a(t) = \frac{1}{2 \sqrt{2}} \ln(t).
$$
Cela permet d'obtenir le développement asymptotique souhaité pour $\med^*(t)$.

\bigskip

Revenons maintenant au cas du mouvement brownien branchant.
Notons que la variable aléatoire 
$$
N_{q(t)}(t) = \sum_{k=1}^{N(t)} 1_{X_k(t) \ge q(t)}
$$
admet le même premier moment que la variable aléatoire $N_{q(t)}^*(t)$.
Le premier moment de $N_{q(t)}(t)$ ne donne par contre pas ici une information précise sur la probabilité $P(X_1(t) \ge q(t))$.
En effet, si une particule parvient au-delà de $q(t)$ à l'instant $t$, 
plusieurs de ses ancêtres sont probablement à une altitude élevée et plusieurs de leurs descendants sont donc probablement au-delà de $q(t)$.
Ainsi, le $q(t)$ pour lequel $E(N_{q(t)}(t))$ est d'ordre $1$ surestime la valeur de la médiane $\med(t)$ de $X_1(t)$.

\`A un instant $s \le t$, l'espérance du nombre de particules du mouvement brownien branchant est $e^s$.
En appliquant brutalement une borne de premier moment on obtient donc que, à l'instant $s$, les particules ne peuvent pas être très loin au-dessus de $\sqrt{2}s$.
L'une des idées est de rajouter ce type de contraintes sur la trajectoire des particules d'intérêt.
Plus précisément, 
appelons particule basse une particule au-delà de $q(t)$ à l'instant $t$ et dont la trajectoire est restée en-dessous de la droite $s\mapsto q(t)st^{-1}+2$.
Appelons particule haute une particule au-delà de $q(t)$ à l'instant $t$ et dont la trajectoire a touché la droite d'équation $s\mapsto q(t)st^{-1}+1$.
Notons $B_{q(t)}(t)$ le nombre de particules basses et $H_{q(t)}(t)$ le nombre de particules hautes.
On a :
$$
\frac{E(B_{q(t)}(t))^2}{E(B_{q(t)}(t)^2)} \le P(B_{q(t)}(t) \ge 1) \le P(N_{q(t)}(t) \ge 1) \le E(B_{q(t)}(t)) + P(H_{q(t)}(t) \ge 1).
$$
La majoration sur la position des ancêtres des particules basses permet de majorer efficacement l'espérance de $B_{q(t)}(t)$ 
conditionnellement à l'existence d'une particule basse.
Cela permet - pour des choix de $q(t)$ pertinents - de contrôler le second moment de $B_q(t)$ par le carré de son premier moment.
Par ailleurs, toute particule haute admet un ancêtre sur la droite d'équation $s \mapsto q(t)st^{-1}+1$.
Cela permet - toujours pour des choix de $q(t)$ pertinents - de minorer l'espérance de $B_{q(t)}(t)$ conditionnellement à l'existence d'une particule haute.
La probabilité d'existence d'une particule haute peut ainsi être contrôlée par l'espérance de $B_{q(t)}(t)$.
La conclusion est que la détermination de la médiane $\med(t)$ de $X_1(t)$ se ramène à la détermination du paramètre $q(t)$ 
pour lequel le premier moment de $B_q(t)$ est d'ordre $1$.

Estimons maintenant $B_{q(t)}(t)$. 
Essentiellement, l'introduction de la contrainte sur la trajectoire divise par $t$ le premier moment :
\begin{equation}\label{e:leading}
E(B_{q(t)}(t)) \sim \frac{C}{t} E(N_{q(t)}(t)) = \frac{C}{t} E(N^*_{q(t)}(t)).
\end{equation}
Le facteur $1/t$ se comprend facilement dans un cadre discret en temps avec des évènements légèrement différents.
Si $X_1, \dots, X_t$ sont des v.a.i.i.d. de loi commune diffuse, alors :
$$
P\left(X_1+ \dots + X_s \le \frac st \; X_t \text{ pour tout } s\le t \Big|X_1 + \cdots + X_t \ge q(t)\right) = \frac1t.
$$
La preuve est une conséquence des deux remarques suivantes : 
les rotations des $X_i$ laissent invariant l'évènement par lequel on conditionne ;
l'autre évènement est vérifié pour une unique rotation des $X_i$.
Rappelons que $N_{q(t)}(t)$ et $N^*_{q(t)}(t)$ ont le même premier moment.
En combinant les relations \eqref{e:m1indep} et \eqref{e:leading}  on obtient, 
toujours lorsque $a(t)$ est négligeable devant $\sqrt{t}$ :
\begin{equation} \nonumber
E \big(B_{\sqrt{2}t - a(t)} (t) \big) \sim  \frac{C}{t^{3/2}} e^{a(t)\sqrt{2}}.
\end{equation}
Le premier moment de $B_q(t)$ est ainsi d'ordre $1$ pour $t$ grand lorsque :
$$
q(t) = \sqrt{2}t - \frac{3}{2 \sqrt{2}} \ln(t).
$$
Cela fournit l'estimation souhaitée pour la médiane $\med(t)$.


\section{Mouvement brownien branchant vu depuis son sommet}

\subsection{Convergence}

On pose :
\begin{equation}\label{e:mt}
m(t) = \sqrt{2}t - \frac{3}{2 \sqrt{2}} \ln(t).
\end{equation}
Nous nous intéressons aux processus suivants :
\begin{itemize}
\item le mouvement brownien branchant :
$$
\cN(t) = \sum_{k=1}^{N(t)} \delta_{X_k(t)}.
$$
\item le mouvement brownien branchant vu depuis $m(t)$ :
$$
\cN_m(t) = \sum_{k=1}^{N(t)} \delta_{X_k(t)-m(t)} = \cN(t)-m(t).
$$
\item le mouvement brownien branchant vu depuis $m(t)+2^{-1/2}\ln(C_wZ)$ où $Z$ est la limite de la martingale définie en $\eqref{e:Z}$ 
et où $C_w$ est définie par \eqref{e:asymptotiquew} :
$$
\cN_Z(t) = \sum_{k=1}^{N(t)} \delta_{X_k(t)-m(t)-\frac1{\sqrt{2}}\ln(C_wZ)} = \cN(t)-\left(m(t)+\frac1{\sqrt{2}}\ln(C_wZ)\right).
$$
\end{itemize}
Si $f$ est une fonction de $\R$ dans $\R$ et  si $\cP$ est un processus ponctuel nous noterons souvent $f(\cP)$ l'intégrale de $f$ contre la mesure $\cP$.
Ainsi et par exemple $f(\cN(t))=\sum_{k=1}^{N(t)} f(X_k(t))$.

\bigskip

Le mouvement brownien vu depuis son sommet converge en loi.
\ABetK ont démontré le résultat suivant.

\begin{theo}[\cite{ABK3}] \label{t:ABKnew}
\begin{itemize}
\item Le processus $\cN_m(t)$ converge en loi vers un processus ponctuel $\cN_m$.
\item Le processus ponctuel $\cN_m$ a la même loi que le processus
$$
\sum_{x \in \cP} x+\frac{1}{\sqrt{2}}\ln(C_wZ) + \cC_x
$$
où $\cP$ est un processus de Poisson ponctuel d'intensité $\sqrt{2} \exp(-\sqrt{2}x) dx$ indépendant de $Z$ et où, 
conditionnellement à $Z$ et $\cP$,  $(\cC_x)_{x \in \cP}$ est une de famille de copies indépendantes d'un processus ponctuel $\cC$.
\item Le processus ponctuel $\cC$ vérifie $\max \cC = 0$. Sa loi est celle du processus $\cDABK$ défini par \eqref{e:DABK}.
\end{itemize}
\end{theo}

\ABBetS ont donné peu après la version plus fine suivante.

\begin{theo}[\cite{ABBS}] \label{t:ABBSnew}
\begin{itemize}
\item Le processus $(\cN_Z(t),Z(t))$ converge en loi vers $(\cN_Z,Z)$ où $\cN_Z$ est un processus ponctuel indépendant de $Z$.
\item Le processus $\cN_Z$ a la même loi que le processus
\begin{equation}\label{e:dimanche}
\sum_{x \in \cP} x+\cC_x.
\end{equation}
où $\cP$ est un processus de Poisson ponctuel d'intensité $\sqrt{2} \exp(-\sqrt{2}x) dx$ et où, 
conditionnellement à  $\cP$,  $(\cC_x)_{x \in \cP}$ est une de famille de copies indépendantes d'un processus ponctuel $\cC$.
\item Le processus ponctuel $\cC$ vérifie $\max \cC = 0$. Sa loi est celle du processus $\cQ$ défini par \eqref{e:defYcQ}.
\end{itemize}
\end{theo}

Des résultats précédents on déduit facilement que le processus limite admet une propriété de superposabilité : 
la superposition de copies indépendantes de $\cN_Z$ suit, à une translation près, la même loi que $\cN_Z$.
Maillard \cite{Maillard} a établi que la superposabilité caractérise les processus de la forme \eqref{e:dimanche}.

Les résultats principaux de \cite{ABK3} et \cite{ABBS} concernent la description du processus limite et, en particulier, les descriptions de la décoration $\cC$
données dans les deux parties suivantes.

\bigskip

Les propriétés du mouvement brownien branchant vu depuis un extrémum ont également été étudiées récemment par Brunet-Derrida 
\cite{Brunet-Derrida-EPL, Brunet-Derrida-JSP}.
Ces travaux reposent notamment sur le lien entre mouvement brownien et équation F-KPP.
On y trouve - au moins implicitement - la preuve de la convergence ainsi que différentes conjectures.
Certaines sont établies, comme la superposabilité de la limite et la caractérisation de Maillard.
D'autres sont encore ouvertes, comme celles pourtant sur la densité du processus limite ou sur l'écart moyen entre deux particules successives du processus limite.

\bigskip

L'approche de \ABetK est de nature plutôt analytique et s'inspire en partie des travaux récents de Brunet et Derrida \cite{Brunet-Derrida-EPL,Brunet-Derrida-JSP}.
L'approche de \ABBetS est de nature plutôt probabiliste et repose en partie sur les résultats récents de \ABetK \cite{ABK1} décrits dans la partie \ref{s:ABK1}.

\bigskip

Dans \cite{ABK5}, \ABetK établissent la convergence en moyenne ergodique de $\cN_m(t)$ vers $\cN_m$.
Les mêmes auteurs avaient établi précédemment dans \cite{ABK4} la convergence en moyenne ergodique de $X_1(t)-m(t)$ vers la somme d'une variable aléatoire
de loi Gumbel et de $2^{-1/2}\ln(C_wZ)$.
La preuve repose notamment sur une extension du résultat de Lalley et Sellke \eqref{e:LS} et sur des résultats sur la généalogie 
du processus permettant de controler certaines corrélations.

\subsection{Loi de la décoration $\cC$ : description de \ABetK}

\ABetK établissent le résultat suivant :

\begin{theo}[\cite{ABK3}]\label{t:ABK-conditionne} Soient $a > 0$ et $b \in \R$.
\begin{itemize}
\item Conditionnellement à 
\begin{equation}\label{e:bonbon}
\{ X_1(t) - (\sqrt{2}t+a\sqrt{t}+b) > 0\},
\end{equation}
le processus 
$$
\sum_{k=1}^{N(t)} \delta_{X_k(t)-(\sqrt{2}t+a\sqrt{t}+b)}
$$
converge en loi vers un processus non trivial $\overline{\cN}$.
\item La loi du processus limite $\overline{\cN}$ ne dépend pas du choix de $a$ et de $b$.
\item Le maximum de $\overline{\cN}$ suit une loi exponentielle de moyenne $2^{-1/2}$.
\end{itemize}
\end{theo}

Rappelons que, sans conditionnemement, $\max_k X_k(t) - \sqrt{2}t $ tend vers $-\infty$.
Posons :
\begin{equation}\label{e:DABK}
\cDABK = \overline{\cN} - \max(\overline{\cN}).
\end{equation}
Ce processus ponctuel est la décoration apparaissant dans le théorème \ref{t:ABKnew}.
Du théorème précédent, on déduit que, conditionnellement à \eqref{e:bonbon}, le processus 
$$
\left(\cN(t)-X_1(t), X_1(t) - (\sqrt{2}t+a\sqrt{t}+b)\right)
$$
converge en loi vers $(\cDABK,H)$ où $\cDABK$ et $H$ sont indépendants et $H$ suit une loi exponentielle de moyenne $2^{-1/2}$.

\subsection{Loi de la décoration $\cC$ : description de \ABBetS}

Pour tout $0 \le s \le t$, notons $X_{1,t}(s)$ la position à l'instant $s$ de la particule qui est en $X_1(t)$ à l'instant $t$.
Autrement dit, $X_{1,t}(\cdot)$ est la trajectoire de la particule qui réalise le maximum à l'instant $t$.
Nous nous intéressons au renversement de cette trajectoire que nous notons $Y_t(\cdot)$.
Il est défini en $s \in [0,t]$ par : 
$$
Y_t(s) = X_{1,t}(t-s) - X_1(t).
$$
Si $1 \le i,j \le N(t)$, nous notons $\tau_{i,j}(t)$ l'instant de $[0,t]$ auquel les particules $X_i(t)$ et $X_j(t)$ se sont séparées.
Pour tout $\zeta$, nous considérons le processus ponctuel suivant :
\begin{equation}\label{e:cQzeta}
\cQ(t,\zeta) = \sum_{1 \le k \le N(t) : \tau_{1,k}(t) > t-\zeta} \delta_{X_k(t)-X_1(t)}.
\end{equation}
C'est, vu depuis $X_1(t)$, l'ensemble des particules à l'instant $t$ qui se sont séparées de $X_1(t)$ après l'instant $t-\zeta$.

\ABBetS établissent la convergence de $(Y_t,\cQ(t,\zeta),X_1(t)-m(t))$ lorsque $t$ puis $\zeta$ tendent vers l'infini.
Décrivons tout d'abord l'objet limite.
Pour cela, nous commençons par définir pour tout $b>0$ un processus à valeurs réels $(\Gamma^{(b)}_t)_{t \ge 0}$.
Soient  $(B(t))_{t \ge 0}$  un mouvement brownien standard et $(R(t))_{t\ge 0}$ un processus de Bessel indépendant de dimension $3$ issu de $0$. 
Notons $T_b$ le premier temps d'atteinte du niveau $b$ par le mouvement brownien $B$.
Le processus $\Gamma^{(b)}$ est défini ainsi :
\begin{equation}
\Gamma^{(b)}(t) := \begin{cases} B(t), &\text{ si $t\in [0, \, T_b]$}, \cr\cr b-R(t-T_b), &\text{ si $t\ge T_b$.} \cr \end{cases}
\end{equation}
Nous définissons alors $(Y^{(b)}(t))_{t \ge 0}$ par $Y^{(b)}(t)=\Gamma^{(b)}-\sqrt{2}t$.
Conditionnellement à $Y^{(b)}$ nous définissons un processus ponctuel $\overline{\cQ}^{Y^{(b)}}$ par 
$$
\overline{\cQ}^{Y^{(b)}} = \delta_0 + \sum_{t \in \chi} \left(\cN_t(t)+Y^{(b)}(t)\right),
$$
où $\chi$ est un processus de Poisson ponctuel indépendant sur $[0,+\infty[$ de mesure d'intensité $2dt$ et où, conditionnellement à ce qui précède,
les $(\cN_t(\cdot))_{t \in \chi}$ sont des copies indépendantes du mouvement brownien branchant.
Nous définissons alors la loi de $(Y,\cQ)$ par :
\begin{equation}\label{e:defYcQ}
E \phi(Y,\cQ) = C\int_0^{+\infty} db E\left( \phi\left(Y^{(b)},\overline{\cQ}^{Y^{(b)}}\right) 1_{\max \overline{\cQ}^{Y^{(b)}} \le 0}\right)
\end{equation}
pour tout $\phi$.
La quantité $C$ est la constante de normalisation.
Le processus ponctuel $\cQ$ est la décoration apparaissant dans le théorème \ref{t:ABBSnew}.

\ABBetS établissent le résultat suivant.

\begin{theo}[\cite{ABBS}]\label{t:ABBS-2.3}
$$
\lim_{\zeta \to \infty} \lim_{t \to \infty} (Y_t , \cQ(t,\zeta), X_1(t)-m(t)) = (Y,\cQ,W) \text{ en loi}
$$
où $(Y,\cQ)$ et $W$ sont indépendants. 
\end{theo}

\section{Localisation des trajectoires des particules extrémales et généalogie}
\label{s:ABK1}

Nous décrivons rapidement dans cette partie les travaux de \ABetK obtenus dans \cite{ABK1}.
Ces résultats sont utilisés par \ABBetS dans leur preuve du théorème \ref{t:ABBSnew}.

\subsection{Trajectoires des particules extrémales}

Soient $t \ge 0$ et $\alpha > 0$.
Les fonctions $U^+_{t,\alpha}$ et $U^-_{t,\alpha}$ sont définies, pour tout $s \in [0,t]$, par :
$$
U^+_{t,\alpha}(s) = \frac st m(t) + \min(s^\alpha, (t-s)^{\alpha})
$$
et
$$
U^-_{t,\alpha}(s) = \frac st m(t) - \min(s^\alpha, (t-s)^{\alpha}).
$$

Rappelons que, pour tout $k \le N(t)$, $X_{k,t}(\cdot)$ est la trajectoire entre les instants $0$ et $t$ de la particule en $X_k(t)$ à l'instant $t$.
\ABetK établissent le résultat suivant.

\begin{theo}[\cite{ABK1} - Théorèmes 2.2, 2.3 et 2.5] \label{t:localisation}
Soit $A>0$.
Soit $0<\alpha < 1/2$.
\begin{enumerate}
\item 
$$
\lim_{r \to \infty} \sup_{t \ge 3r} P\big( \exists  k \le N(t), \exists s \in [r,t-r] : X_{k,t}(s)  \ge U^+_{t,\alpha}(s) \big) = 0.
$$
\item 
$$
\lim_{r \to \infty} \sup_{t \ge 3r} 
P\big( \exists k \le N(t) : |X_k(t)-m(t)| \le A \text{ et } \exists s \in [r,t-r] : X_{k,t}(s)  \ge U^-_{t,1/2-\alpha}(s) \big) = 0.
$$
\item 
$$
\lim_{r \to \infty} \sup_{t \ge 3r} 
P\big( \exists k \le N(t) : |X_k(t)-m(t)| \le A \text{ et } \exists s \in [r,t-r] : X_{k,t}(s)  \le U^-_{t,1/2+\alpha}(s) \big) = 0.
$$
\end{enumerate}
\end{theo}

Notons que le premier point concerne toutes les trajectoires.
Il pourrait par ailleurs être également énoncé ainsi :
$$
\lim_{r \to \infty} \sup_{t \ge 3r} P\big(\exists s \in [r,t-r] : X_1(s)  \ge U^+_{t,\alpha}(s) \big) = 0.
$$
C'est un contrôle uniforme en temps sur le maximum du mouvement brownien branchant.
La preuve du premier point donnée par \ABetK repose sur le résultat suivant de Bramson \cite{Bramson1}.
Il existe une constante $C$ telle que, pour tout $0<y<\sqrt{t}$ :
$$
P(X_1(t) \ge m(t) + y) \le C(1+y)^2\exp(-\sqrt{2}y).
$$
Cette inégalité permet d'obtenir simplement un contrôle uniforme sur les temps entiers.
Un contrôle sur tous les temps s'en déduit par des majorations grossières.

\bigskip

Les deux derniers points ne concernent que les trajectoires des particules proches de $m(t)$ à l'instant $t$.
Ils se déduisent relativement facilement du premier point qui affirme que ces particules suivent essentiellement la trajectoire
d'une excursion brownienne en-dessous de $s \mapsto (s/t)m_t$.
Notons qu'il existe par contre, avec grande probabilité, des particules proches de $m(t/2) = \frac12m(t) + O(\ln(t))$ à l'instant $t/2$
et proches de $2m(t/2) = m(t) + O(\ln(t))$ à l'instant $t$. 

La forme de la trajectoire des particules extrémales peut également se comprendre comme le résultat d'une compétition entre energie 
- liée à la probabilité qu'une particule à une hauteur $q$ à un instant $s$ fournisse un descendant extrémal à l'instant $t$ -
et entropie -liée au nombre de particules à une hauteur $q$ à un instant $s$.

\ABBetS établissent dans \cite{ABBS} une version légèrement différente du théorème \ref{t:localisation}.

\subsection{Généalogie du mouvement brownien branchant}

Rappelons que $\tau_{i,j}(t)$ désigne l'instant auquelle les particules $X_i(t)$ et $X_j(t)$ se sont séparées.
Dans \cite{ABK1}, \ABetK établissent le résultat suivant.

\begin{theo}[\cite{ABK1} - Théorème 2.1] \label{t:genealogie}
Soit $A>0$. On a :
$$
\lim_{r \to \infty} \sup_{t \ge 3r} P\big(\exists i,j \le N(t) : X_i(t),X_j(t) \in [m(t)-A,m(t)+A] \hbox{ et } \tau_{i,j}(t) \in [r,t-r]\big)=0.
$$
\end{theo}
La preuve de ce théorème repose sur les résultats de localisation des trajectoires de particules extrémales énoncés dans le théorème \ref{t:localisation}.
Observons les particules à un instant $s$ loin de $0$ et de $t$.
Les particules susceptibles de fournir un descendant extrémal à l'instant $t$ sont situées aux environs de $st^{-1}m(t) - \min(s,(t-s))^{1/2}$.
La probabilité que l'une d'entre elle fournisse un descendant extrémal à l'instant $t$ est faible.
Le nombre de telles particules est par contre élevé.
Ces deux facteurs se compensent.
Par contre, 
la probabilité que l'une de ces particules se scinde en deux et que chacune de ces deux particules admette un descendant extrémal à l'instant $t$ est faible.

Techniquement, la preuve repose sur une majoration de l'espérance du nombre de couples $(i,j)$ vérifiant les conditions apparaissant dans l'énoncé du théorème
ainsi que les conditions de localisation des trajectoires.

\section{Preuves}

\subsection{Approches de \ABetK}

\subsubsection{\'Equations F-KPP et mouvement brownien branchant}
Notons pour commencer que si une application $v$ est solution de l'équation F-KPP
\begin{equation}\label{e:FKPPv}
\frac{\partial v}{\partial t}
=
\frac12 \frac{\partial^2 v}{\partial^2 x} + v^2-v
\end{equation}
alors $u=1-v$ est solution de l'équation F-KPP
\begin{equation}\label{e:FKPPu}
\frac{\partial u}{\partial t}
=
\frac12 \frac{\partial^2 u}{\partial^2 x} + u-u^2.
\end{equation}
Voici maintenant une version plus précise de l'observation - élémentaire et féconde - de McKean liant le mouvement brownien branchant et les équations F-KPP.

\begin{lemm}[\cite{McKean}] \label{l:McKean} Soit $f:\R\to[0,1]$. La fonction $v$ définie par
$$
v(t,x) = E\left(\prod_{k=1}^{N(t)}f(x+X_k(t))\right)
$$
est solution de l'équation F-KPP \eqref{e:FKPPv} de condition initiale $v(0,\cdot)=f$.
\end{lemm}

Explicitons $v$ dans les cas qui nous intéresseront principalement dans la suite.
Soit $\phi:\R\to\R$ une fonction positive, continue et à support compact.
Soit $\delta$ un nombre réel.
Prenons pour condition initiale la fonction $f$ définie par :
$$
f(y) = \exp(-\phi(-y))1_{]-\infty,\delta]}(-y).
$$
La fonction $v$ associée vérifie alors :
\begin{eqnarray*}
v(t,x)  
 & = & E\left(\prod_{k=1}^{N(t)}\exp(-\phi(-x-X_k(t))1_{]-\infty,\delta]}(-x-X_k(t))\right), \\
 & = & E\left(\prod_{k=1}^{N(t)}\exp(-\phi(-x+X_k(t))1_{]-\infty,\delta]}(-x+X_k(t))\right), \\
 & = & E\left(\exp(-\phi(\cN(t)-x)1_{X_1(t) \le \delta+x}\right).
\end{eqnarray*}

Si $\phi$ est constante égale à $0$ et si $\delta=0$ on a par exemple :
\begin{equation} \label{e:vmax}
v(t,m(t)+x) = P(X_1(t) \le m(t)+x).
\end{equation}
L'étude de la convergence en loi de $X_1(t)-m(t)$ se ramène ainsi à celle de la convergence des fonctions $t \mapsto v(t,m(t)+x)$ pour tout réel $x$.

Si l'on considère la condition initiale $f$ définie par $f(y) = \exp(-\phi(-y))$ alors on a par exemple :
$$
v(t,m(t)) = E(\exp(-\phi(\cN_m(t)))).
$$
\'Etudier la convergence en loi du processus $(\cN_m(t))_t$ revient ainsi à étudier la convergence des fonctions $t \mapsto v(t,m(t))$ 
associées aux différentes fonctions $\phi$.



\subsubsection{Condition initale : front borné}

Dans cette partie, nous considérons le comportement asymptotique des solutions de l'équation F-KPP \eqref{e:FKPPu} de condition initiale $g:\R\to[0,1]$ 
satisfaisant :
\begin{equation}\label{e:condition}
\hbox{Il existe }A>0\hbox{ tel que }g(x)=1\hbox{ pour }x \le -A \hbox{ et }g(x)=0\hbox{ pour } x \ge A.
\end{equation}

Bramson a établi dans \cite{Bramson2} le résultat de convergence suivant, dans lequel $m(t)$ est défini par \eqref{e:mt}
et où $w$ est la fonction de répartition de la variable aléatoire $W$ définie en \eqref{e:W}. 

\begin{theo}[\cite{Bramson2}] \label{t:CSCVU}
Soit $u$ une solution de l'équation F-KPP \eqref{e:FKPPu} de condition initiale $u(0,\cdot)=g$ où 
$g:\R\to[0,1]$ est une fonction satisfaisant la condition \eqref{e:condition}.
Alors il existe une constante $C'(g)$ telle que
$$
u(t,x + m(t)) \to 1-w(x+C'(g)) \hbox{ uniformément en }x\hbox{ lorsque }t \to \infty.
$$ 
\end{theo}

Par \eqref{e:representationw} on a :
$$
1-w(x+C'(g)) = 1-E\left(\exp\left(-C''(g)Ze^{-\sqrt{2}x}\right)\right)
$$
où $C''(g)=C_w e^{-\sqrt{2}C'(g)}$.
La constante $C''(g)$ peut se retrouver à partir du comportement asymptotique de la limite. 
En effet, par \eqref{e:asymptotiquew}, on a :
\begin{equation} \label{e:tgv}
1-w(x+C'(g)) \sim C''(g)xe^{-\sqrt{2}x}.
\end{equation}
Nous noterons souvent ces constantes $C'(u)$ et $C''(u)$ où $u$ est la solution de l'équation F-KPP \eqref{e:FKPPu} de condition initiale $g$.
 
Le comportement asymptotique \eqref{e:bramson-mediane} est une conséquence immédiate des résultats précédents.
En effet, par \eqref{e:vmax} et par le théorème \ref{t:CSCVU} 
(la condition initiale est définie par $g(y)=1-1_{]-\infty,0[}(-y)=1_{]-\infty,0]}$) on obtient :
$$
P(X_1(t) \le m(t)+x) \to w(x+C'(1_{]-\infty,0]})).
$$
Ainsi, $X_1(t)-m(t)-C'(1_{]-\infty,0]})$ converge en loi vers $W$ (on a en fait $C'(1_{]-\infty,0]})=0$).

Bramson fournit dans \cite{Bramson2} une version plus précise du théorème \ref{t:CSCVU}.
Il donne une condition nécessaire et suffisante sur la condition initiale $g$ pour que $u(t,\cdot)$ convenablement translatée converge vers $w$.

\bigskip

Le résultat suivant est une conséquence relativement simple d'un résultat de Bramson également obtenu dans \cite{Bramson2}.
\begin{prop} 
Soit $u$ une solution de l'équation F-KPP \eqref{e:FKPPu} de condition initiale $u(0,\cdot)$ satisfaisant la condition \eqref{e:condition}.
Posons, pour tout $t>r>0$ et tout $x$ réel :
$$
\psi(r,t,x)=\frac{e^{-\sqrt{2}(x-\sqrt{2}t)}}{\sqrt{2\pi(t-r)}}\int_0^{\infty} u(r,y+\sqrt{2}r)e^{y\sqrt{2} - \frac{(y-(x-\sqrt{2}t))^2}{2(t-r)}}
\left(1-e^{-2y\frac{x-m(t)}{t-r}}\right)dy.
$$
Alors, pour tout $r$ suffisamment grand, pour tout $t \ge 8r$ et tout $x \ge m(t)+8r$, on a :
\begin{equation}\label{e:encadrementpsi}
\gamma^{-1}(r) \psi(r,t,x) \le u(t,x) \le \gamma(r)\psi(r,t,x)
\end{equation}
où $\gamma(r)$ décroît vers $1$ lorsque $r$ tend vers l'infini.
\end{prop}

Notons que $\psi(r,\cdot,\cdot)$ ne dépend de $u$ qu'à travers $u(r,\cdot)$.
La proposition donne le comportement de $u(t,x)$ pour $t$ grand et $x-m(t)$ grand.
Explicitons ce comportement dans différents régimes.

Posons :
\begin{equation}\label{e:Cr}
C(r,u) = \sqrt{\frac{2}{\pi}} \int_0^{\infty} u(r,y+r\sqrt{2})ye^{y\sqrt{2}}dy.
\end{equation}
\ABetK vérifient que cette quantité est finie par comparaison de $u$ avec la solution explicite de l'équation F-KPP \eqref{e:FKPPu} linéarisée. 
On en déduit alors simplement, pour tout réel $x$ :
$$
\lim_{t \to \infty} e^{x\sqrt{2}}\frac{t^{3/2}}{\frac{3}{2\sqrt{2}}\ln(t)} \psi(r,t,\sqrt{2}t+x) = C(r,u).
$$
L'encadrement \eqref{e:encadrementpsi} et le fait que $\gamma$ tende vers $1$ permet alors d'en déduire la convergence de $C(r,u)$ vers un réel 
strictement positif $C(u)$,
\begin{equation}\label{e:definitionCu}
C(u) = \lim_{r \to \infty} \sqrt{\frac{2}{\pi}} \int_0^{\infty} u(r,y+r\sqrt{2})ye^{y\sqrt{2}}dy,
\end{equation}
ainsi que la convergence suivante, pour tout réel $x$ :
$$
\lim_{t \to \infty} e^{x\sqrt{2}}\frac{t^{3/2}}{\frac{3}{2\sqrt{2}}\ln(t)} u(t,\sqrt{2}t+x) = C(u).
$$
Les arguments précédents proviennent de Chauvin-Rouault \cite{Chauvin-Rouault-90}.

On montre similairement
$$
\lim_{t \to \infty} \frac{e^{x\sqrt{2}}}{x} \psi(r,t,m(t)+x) = C(r,u),
$$
puis
$$
\lim_{x\to\infty} \lim_{t \to \infty} \frac{e^{x\sqrt{2}}}{x} u(t,m(t)+x) = C(u).
$$
Par \eqref{e:tgv}, on en déduit $C(u)=C''(u)$.
En utilisant \eqref{e:representationw} nous obtenons donc la représentation intégrale suivante :
$$
\lim_{t\to\infty} u(t,m(t)+x) = 1-E\left(\exp\left(-C(u)Ze^{-\sqrt{2}x}\right)\right)
$$
où $C(u)$ est défini par \eqref{e:definitionCu}.

Nous aurons également besoin du comportement asymptotique des solutions dans le régime suivant.
Soient $a>0$ et $b$ un nombre réel.
On obtient comme précédemment le résultat suivant :
\begin{equation} \label{e:asymptotique-racine}
\lim_{t \to \infty} \frac{e^{\sqrt{2}(a\sqrt{t}+b)}t^{3/2}}{a\sqrt{t}+b} \psi(r,t,\sqrt{2}t+a\sqrt{t}+b) = C(r,u)e^{-a^2/2}.
\end{equation}

\bigskip

Soit $A$ un réel tel que $u(0,\cdot) \le 1_{]-\infty,A]}$. En utilisant le lemme \ref{l:McKean} on obtient :
$$
u(t,x) \le P(X_1(t) \ge x-A).
$$
Un contrôle suffisamment fin sur la queue de la distribution du maximum permet d'en déduire :
$$
\lim_{A_1 \searrow 0} \limsup_{r \to \infty} \int_0^{A_1\sqrt{r}} u(r,y+\sqrt{2}r) ye^{y\sqrt{2}}dy=0
$$
et
$$
\lim_{A_2 \nearrow \infty} \limsup_{r \to \infty} \int_{A_2\sqrt{r}}^{\infty} u(r,y+\sqrt{2}r) ye^{y\sqrt{2}}dy=0.
$$
Résumons ces propriétés avec l'écriture suivante :
\begin{equation} \label{e:racinev1}
C(u) = \lim_{r \to \infty} \int_{y \approx \sqrt{r}}  u(r,y+\sqrt{2}r) ye^{y\sqrt{2}}dy.
\end{equation}

Ces contrôles permettent notamment d'obtenir simplement le comportement de $C(u)$ lorsque la condition initiale est translatée :
\begin{equation} \label{e:Ctranslate}
C(u(\cdot,\cdot+x)) = C(u)e^{-\sqrt{2}x}.
\end{equation}


\subsubsection{Condition initiale localisée}

Certains des résultats précédents restent valides pour des conditions initiales différentes.

\begin{prop} \label{p:3.2} Soit $\phi:\R\to\R$ une fonction positive, continue, à support compact et non identiquement nulle.
Notons $u$ la solution de l'équation F-KPP \eqref{e:FKPPu} associée à la condition initiale $g=1-\exp(-\phi(-\cdot))$.
On a :
\begin{equation} \label{e:limite}
\lim_{t \to \infty} E\left(e^{-\phi(\cN_m(t))}\right) = E\left(e^{-C(u)Z}\right)
\end{equation}
où
\begin{equation}\label{e:Cugeneral}
C(u)=\lim_{t\to\infty} \sqrt{\frac{2}{\pi}} \int_0^{\infty} u(t,t\sqrt{2}+y)ye^{y\sqrt{2}}dy
\end{equation}
est une constante strictement positive ne dépendant que de $\phi$.
Par ailleurs, pour tout réel $x$ on a :
\begin{equation}\label{e:Ctranslategeneral}
C(u(\cdot,\cdot+x)) = C(u)e^{-\sqrt{2}x}.
\end{equation}
On a enfin, avec le même sens qu'en \eqref{e:racinev1} : 
\begin{equation} \label{e:racinegeneral}
C(u) = \lim_{r \to \infty} \int_{y \approx \sqrt{r}}  u(r,y+\sqrt{2}r) ye^{y\sqrt{2}}dy.
\end{equation}
\end{prop}

Cette proposition permet d'obtenir la convergence en loi du processus $(\cN_m(t))_t$. 
Elle donne également des informations sur le processus limite.
Les relations \eqref{e:limite} à \eqref{e:Ctranslategeneral} s'inspirent des travaux de Chauvin-Rouault \cite{Chauvin-Rouault-90}.

On a pour tout nombre réel $x$ :
$$
\lim_{t \to \infty} u(t,m(t)+x) = \lim_{t \to \infty} 1-E\left(e^{-\phi(\cN_m(t)-x)}\right) = 1-E\left(e^{-C(u)Ze^{-\sqrt{2}x}}\right).
$$
La convergence n'est pas uniforme en $x$. Par symétrie, on observe en effet la formation d'un front symétrique autour de la position $-m(t)$.

La limite précédente peut s'écrire :
$$
E\left(e^{-Ze^{-\sqrt{2}\left(x-\frac1{\sqrt{2}}\ln(C(u))\right)}}\right).
$$
Le terme $\frac1{\sqrt{2}}\ln(C(u)$ s'interprète ainsi comme un terme de retard de l'onde stationnaire.

La limite apparaissant dans \eqref{e:limite} et la relation \eqref{e:Cugeneral} permettent de donner une description probabiliste du processus limite.
C'est l'objet du théorème \ref{t:auxiliaire}. Cette description fait intervenir le processus auxiliaire de la partie suivante. 
La relation \eqref{e:racinegeneral} permet de donner une description probabiliste alternative du processus limite.
C'est la description qui fait l'objet du théorème \ref{t:ABKnew}.

\bigskip
Dans la suite de cette partie, nous donnons quelques détails sur la preuve de la proposition \ref{p:3.2}.
Pour tout réel $\delta$, on considère la condition initiale $g_\delta$ définie par $g_\delta(x)=1-\exp(-\phi(-x))1_{-x \le \delta}$.
La condition initiale $g_\delta$ vérifie la condition \eqref{e:condition}.
Notons $u_\delta$ la solution de l'équation \eqref{e:FKPPu} associée.
La preuve est naturelle et consiste à vérifier que $u_\delta$ converge convenablement vers $u$.

Le lemme \ref{l:McKean} permet d'obtenir facilement l'encadrement suivant :
\begin{equation}\label{e:concert}
u_\delta(t,x) -  P(X_1(t) > \delta+x) \le u(t,x) \le u_\delta(t,x).
\end{equation}
Ainsi :
$$
C(t,u_\delta) - C(t,1_{]-\infty,-\delta]}) \le C(t,u) \le C(t,u_\delta).
$$
En prenant la limite en $t$ et en $\delta$ et en utilisant \eqref{e:Ctranslate} on obtient :
$$
\lim_{t \to \infty} C(t,u) = \lim_{\delta \to \infty} \lim_{t \to \infty} C(t,u_\delta),
$$
c'est-à-dire, en appelant $C(u)$ la première limite :
$$
C(u) = \lim_{\delta\to\infty} C(u_\delta).
$$
Notons que $\delta \mapsto C(u_\delta)$ est croissante.
La relation \eqref{e:Ctranslategeneral} se déduit alors de \eqref{e:Ctranslate}.
L'encadrement \eqref{e:concert} permet également d'obtenir :
$$
u_\delta(t,m(t)+x) - P(X_1(t) > m(t)+x+\delta) \le u(t,m(t)+x) \le u_\delta(t,m(t)+x).
$$
En prenant la limite en $t$ puis en $\delta$ dans l'encadrement, on obtient :
$$
\lim_t u(t,m(t)+x) = 1-\lim_\delta E\left(\exp\left(-C(u_\delta)Ze^{-\sqrt{2}x}\right)\right) = 1-E\left(\exp\left(-C(u)Ze^{-\sqrt{2}x}\right)\right).
$$
Si $C(u)$ était nul, alors la limite précédente serait nulle pour tout nombre réel $x$.
Or :
$$
u(t,m(t)+x) = 1-E\left(e^{-\phi(\cN_m(t)-x)}\right) \ge 1-E\left(e^{-\phi(X_1(t)-m(t)-x)}\right).
$$
La convergence en loi de $X_1(t)-m(t)$ permet de conclure que $C(u)$ est non nul.

\subsubsection{Processus auxiliaire}

\ABetK introduisent un processus auxiliaire $(\Pi(t))_t$ et montrent qu'il admet la même limite que $(\cN_m(t))_t$.
Ce processus est défini de la manière suivante.
Soit $\eta$ un processus de Poisson ponctuel sur $\R$ de mesure d'intensité 
$$
\sqrt{\frac2\pi}(-x)e^{-\sqrt{2}x}1_{]-\infty,0[}(x)dx.
$$
Conditionnellement à $\eta$, on se donne une famille $(\cNlineaire^x)_{x \in \eta}$ de copies indépendantes du processus $(\cNlineaire(t))_t$ défini par 
$$
\cNlineaire(t)=\cN(t)-\sqrt{2}t.
$$
Le processus ponctuel $\eta$ et la famille $(\cNlineaire^x)_{x \in \eta}$ sont indépendants de $Z$.
Le processus auxiliaire est défini par :
$$
\Pi(t) = \sum_{x \in \eta} \frac1{\sqrt{2}}\ln(Z)+x+\cNlineaire^x(t).
$$

Rappelons que le maximum du processus $\cNlineaire$ converge presque sûrement vers $-\infty$.
Le processus $\Pi(t )$ converge néanmoins vers un processus ponctuel non trivial.
Plus précisément, \ABetK établissent le résultat suivant :

\begin{theo} \label{t:auxiliaire} Les processus $(\Pi(t))_t$ et $(\cN_m(t))_t$ convergent en loi vers la même limite.
\end{theo}

L'étude du comportement asymptotique de $\cN_m(t)$ {\em via} l'introduction de ce processus auxiliaire se rapproche
par certains aspects de la {\em cavity method} introduite par Parisi et ses coauteurs \cite{Mezard-al}.
Nous renvoyons à \cite{ABK3} pour plus de détails.

\bigskip

Montrons ce résultat. Soit $\phi:\R\to\R$ une fonction continue, positive et à support compact.
Notons $u$ la solution de l'équation F-KPP \eqref{e:FKPPu} associée à la condition initiale $1-\exp(-\phi(-\cdot))$.
Rappelons que l'on a :
$$
u(t,y) = 1-E\left(e^{-\phi(-y+\cN(t))}\right)
$$
et donc :
$$
u\left(t,-x + \sqrt{2}t - \frac1{\sqrt{2}}\ln Z\right) = 1-E\left(e^{-\phi(x+\frac1{\sqrt{2}}\ln Z+\cNlineaire(t))}\left|Z\right.\right).
$$
En utilisant ce qui précède, des conditionnement successifs et en utilisant la forme de la transformée de Laplace d'un processus de Poisson ponctuel on obtient :
\begin{eqnarray*}
E\left(e^{-\phi(\Pi(t))}\right) 
 & = & E\left(\exp\left(-\int_{-\infty}^0 u\left(t,-x + \sqrt{2}t - \frac1{\sqrt{2}}\ln Z\right) \sqrt{\frac2\pi}(-x)e^{-\sqrt{2}x}dx\right)\right), \\
 & = & E\left(\exp\left(-\int_0^\infty u\left(t,x + \sqrt{2}t -\frac1{\sqrt{2}}\ln Z\right) \sqrt{\frac2\pi}xe^{\sqrt{2}x}dx\right)\right).
\end{eqnarray*}
En utilisant la proposition \ref{p:3.2} on obtient alors :
\begin{eqnarray*}
\lim_{t\to\infty} E\left(e^{-\phi(\Pi(t))}\right), 
 & = & E\left(\exp\left(-C\left(u\left(\cdot,\cdot -\frac1{\sqrt{2}}\ln Z\right)\right)\right)\right), \\
 & = & E\left(\exp\left(-ZC(u)\right)\right), \\
 & = &  \lim_{t\to\infty} E\left(e^{-\phi(\cN_m(t))}\right) 
\end{eqnarray*}
ce qui conclut la preuve.

\subsubsection{Preuve du théorème \ref{t:ABK-conditionne} et des deux derniers points du théorème \ref{t:ABKnew}}

Commençons par donner les grandes lignes de la preuve du théorème \ref{t:ABK-conditionne}.
Soient $a >0$ et $b \in \R$.
Soit $\phi:\R\to\R$ une fonction continue, positive et à support compact.
On s'intéresse au comportement en grand temps de :
$$
E\left(e^{-\phi(\cN(t)-(\sqrt{2}t+a\sqrt{t}+b))} \left| X_1(t)-(\sqrt{2}t+a\sqrt{t}+b) >0\right.\right).
$$
Pour $\delta>0$, on commence par étudier :
$$
L(t)=E\left(e^{-\phi(\cN(t)-(\sqrt{2}t+a\sqrt{t}+b))} 1_{]-\infty,\delta]}(X_1(t)-(\sqrt{2}t+a\sqrt{t}+b)\left| X_1(t)-(\sqrt{2}t+a\sqrt{t}+b) >0\right.\right).
$$
On peut écrire 
\begin{eqnarray*}
L(t)
 & = & \frac{1-u_2(t,\sqrt{2}t+a\sqrt{t}+b)}{u_3(t,\sqrt{2}t+a\sqrt{t}+b)}-\frac{1-u_1(t,\sqrt{2}t+a\sqrt{t}+b)}{u_3(t,\sqrt{2}t+a\sqrt{t}+b)} \\
 & = & \frac{u_1(t,\sqrt{2}t+a\sqrt{t}+b)}{u_3(t,\sqrt{2}t+a\sqrt{t}+b)}-\frac{u_2(t,\sqrt{2}t+a\sqrt{t}+b)}{u_3(t,\sqrt{2}t+a\sqrt{t}+b)}
\end{eqnarray*}
où $u_1, u_2$ et $u_3$ sont les solutions de l'équation F-KPP \eqref{e:FKPPu} associées aux conditions initiales $f_1, f_2$ et $f_3$ définies par
$f_1(y)=1-\exp(-\phi(-y))1_{]-\infty,0]}(-y)$, $f_2(y)=1-\exp(-\phi(-y))1_{]-\infty,\delta]}(-y)$ et $f_3(y)=1-1_{]-\infty,0]}(-y)$.
En utilisant \eqref{e:encadrementpsi} et \eqref{e:asymptotique-racine} on en déduit :
$$
\lim_{t \to \infty} L(t)= \frac{C(f_1)-C(f_2)}{C(f_3)}.
$$
Notons que l'on ne peut avoir $C(f_1)=C(f_2)$, car sinon on aurait :
$$
\lim_{t\to\infty} u_1(t,m(t)) = \lim_{t\to\infty} u_2(t,m(t))
$$
et donc, en faisant la différence :
$$
\lim_{t\to\infty} E\left(e^{-\phi(\cN(t)-m(t))} 1_{]0,\delta]}(X_1(t)-m(t))\right)=0,
$$
ce qui est exclu compte-tenu du comportement du maximum et de $\cN(t)-m(t)$.
Ainsi, $L(t)$ converge vers un réel strictement positif qui ne dépend pas du choix de $a$ et de $b$.

La loi du maximum de $\cN(t)-(\sqrt{2}t+a\sqrt{t}+b)$ conditionnellement à la stricte positivité de ce maximum s'étudie de manière similaire,
ce qui permet de conclure l'argument de troncature et la preuve du théorème \ref{t:ABK-conditionne}

\bigskip

Dans la fin de cette partie, nous démontrons les deux derniers points du théorème \ref{t:ABKnew}.
Soit $\phi:\R\to\R$ une fonction continue, positive et à support compact.
Notons $u$ la solution de l'équation F-KPP \eqref{e:FKPPu} associée à la condition initiale $1-\exp(-\cdot)$.
Par la proposition \ref{p:3.2} on a :
\begin{equation} \label{e:enfantspostes}
\lim_{t\to\infty} E\left(e^{-\phi(\cN_m(t))}\right) = E\left(e^{-ZC(u)}\right)
\end{equation}
où 
$$
C(u) = \lim_{t\to\infty} \sqrt{\frac{2}{\pi}} \int_0^{\infty} u(t,t\sqrt{2}+y)ye^{y\sqrt{2}}dy.
$$
En utilisant \eqref{e:racinegeneral} on obtient :
\begin{eqnarray*}
C(u) 
 & = & \lim_{t\to\infty} \sqrt{\frac{2}{\pi}} \int_{y \approx \sqrt{t}} u(t,t\sqrt{2}+y)ye^{y\sqrt{2}}dy, \\
 & = & \lim_{t \to \infty} \sqrt{\frac 2\pi} \int_{y \approx \sqrt{t}} E\left(1-\exp\left(-\phi(\cNlineaire(t)-y\right)\right))ye^{\sqrt{2}y}dy.
\end{eqnarray*}
Soit $B$ un réel tel que $\phi \le 1_{[B,+\infty[}$.
On a alors $1-\exp(-\phi(\cNlineaire(t)-y))=0$ si $y < B$.
En écrivant
$$
E\left(1-e^{-\phi(\cNlineaire(t)-y)}\right)
=
E\left(1-e^{-\phi(\cNlineaire(t)-y)} \left| X_1(t) -\sqrt{2}t - y \ge B \right.\right)P(X_1(t) -\sqrt{2}t - y \ge B)
$$
et en utilisant le résultat énoncé à la suite du théorème \ref{t:ABK-conditionne} on obtient :
\begin{eqnarray*}
C(u) 
 & = & \lim_{t \to \infty} 
\int_{y \approx \sqrt{t}} E\left(1-\exp\left(-\phi\left(\cDABK+H+B\right)\right)\right)P(X_1(t) -\sqrt{2}t - y \ge B)\sqrt{\frac 2\pi}ye^{\sqrt{2}y}dy, \\
 & = & E\left(1-\exp\left(-\phi\left(\cDABK+H+B\right)\right)\right)
\lim_{t \to \infty} \int_{y \approx \sqrt{t}} P(X_1(t) -\sqrt{2}t - y \ge B)\sqrt{\frac 2\pi}ye^{\sqrt{2}y}dy. 
\end{eqnarray*}
En utilisant \eqref{e:racinev1} pour la solution de l'équation F-KPP \eqref{e:FKPPu} associée à la condition initiale $1_{]-\infty,-B[}$ on obtient :
\begin{eqnarray*}
C(u)
 & = & E\left(1-\exp\left(-\phi\left(\cDABK+H+B\right)\right)\right)
\lim_{t \to \infty} \int_0^{\infty} P(X_1(t) -\sqrt{2}t - y \ge B)\sqrt{\frac 2\pi} ye^{\sqrt{2}y}dy \\
 & = & E\left(1-\exp\left(-\phi\left(\cDABK+H+B\right)\right)\right) C(1_{]-\infty,B[}).
\end{eqnarray*}
En utilisant \eqref{e:Ctranslate} et le fait que $H$ suit une loi exponentielle de moyenne $2^{-1/2}$, on obtient finalement :
$$
 C(u) = C(1_{]-\infty,0[}) \int_{\R} \left(1-\exp\left(-\phi\left(\cDABK+y\right)\right) \right)\sqrt{2}e^{-\sqrt{2}y}.
$$
De \eqref{e:enfantspostes} et de ce qui précède on déduit les deux derniers points du théorème \ref{t:ABKnew}.


\subsection{Approche de \ABBetS}

\subsubsection{Normalisation}

Pour exposer l'approche de \ABBetS nous préférons reprendre leur normalisation.
Dans la suite, le mouvement brownien a une variance $\sigma^2=2$, une dérive $\rho=2$ et le temps d'atteinte avant fission est toujours de moyenne $1$.
Par ailleurs, nous nous intéressons au minimum et au processus vu depuis son minimum.
Ces normalisations sont telles que :
$$
E\left(\sum_{k=1}^{N(t)} X_k(t)\right)=1 \hbox{ et } E\left(\sum_{k=1}^{N(t)} X_k(t)e^{-X_k(t)}\right)=0.
$$

Nous conservons le nom des objets introduits précédemment, même si les définitions diffèrent parfois légèrement.
Ainsi, nous notons $X_1(t) \le X_2(t) \le \dots \le X_{N(t)}$ les positions des particules classées par ordre {\em croissant}.
On dispose de la convergence en loi de $X_1(t)-m(t)$ vers une variable aléatoire $W$ où
$$
m(t)=\frac32\ln(t).
$$
La fonction de répartition de la variable aléatoire $W$ vérifie :
\begin{equation} \label{e:asymptotiquewABBS}
P(W \le x) \sim_{x \to -\infty} -C_wxe^x.
\end{equation}
La martingale $(Z(t))_t$ est définie par :
$$
Z(t) = \sum_{k=1}^{N(t)} X_i(t)e^{-X_i(t)}.
$$
La martingale converge presque sûrement vers une variable aléatoire strictement positive $Z$.
Le processus $(\cN(t))_t$ est défini comme précédemment.
Le processus $(\cN_Z(t))_t$ est défini par :
$$
\cN_Z(t) = \cN(t)-\left(m(t)-\ln(C_wZ)\right).
$$

On note $X_{1,t}(\cdot)$ la trajectoire de la particule qui réalise le minimum à l'instant $t$.
Le renversement $Y_t(\cdot)$ de cette trajectoire, les instants de fissions $\tau_{i,j}(t)$, le processus ponctuel $\cQ(t,\zeta)$ 
et la famille de processus $(\Gamma^{(b)})_b$ sont également définis comme précédemment.

On pose $Y^{(b)}(t)=-\sigma \Gamma^{(b)}$.
Conditionnellement à $Y^{(b)}$ nous définissons un processus ponctuel $\overline{\cQ}^{Y^{(b)}}$ par 
$$
\overline{\cQ}^{Y^{(b)}} = \delta_0 + \sum_{t \in \chi} \left(\cN_t(t)+Y^{(b)}(t)\right)
$$
où $\chi$ est un processus de Poisson ponctuel indépendant sur $[0,+\infty[$ de mesure d'intensité $2dt$ et où, conditionnellement à ce qui précède,
les $(\cN_t(\cdot))_{t \in \chi}$ sont des copies indépendantes du mouvement brownien branchant.
Nous définissons alors la loi de $(Y,\cQ)$ par :
\begin{equation}\label{e:defYcQ-normalise}
E \phi(Y,\cQ) = C\int_0^{+\infty} db E\left( \phi\left(Y^{(b)},\overline{\cQ}^{Y^{(b)}}\right) 1_{\min \overline{\cQ}^{Y^{(b)}} \ge 0}\right)
\end{equation}
pour tout $\phi$ où $C$ est la constante de normalisation.

\ABBetS établissent les résultats suivants.

\begin{theo}[\cite{ABBS}] \label{t:ABBS-2.3-normalise}
$$
\lim_{\zeta \to \infty} \lim_{t \to \infty} (Y_t , \cQ(t,\zeta), X_1(t)-m(t)) = (Y,\cQ,W) \text{ en loi}
$$
où $(Y,\cQ)$ et $W$ sont indépendants.
\end{theo}

\begin{theo}[\cite{ABBS}] \label{t:ABBS-normalise}
\begin{itemize}
\item Le processus $(\cN_Z(t),Z(t))$ converge en loi vers $(\cN_Z,Z)$ où $\cN_Z$ est un processus ponctuel indépendant de $Z$.
\item Le processus $\cN_Z$ a la même loi que le processus
$$
\sum_{x \in \cP_Z} x+\cC_x
$$
où $\cP_Z$ est un processus de Poisson ponctuel d'intensité $\exp(x) dx$ et où, conditionnellement à $\cP_Z$, $(\cC_x)_{x \in \cP_Z}$
est une suite de copies indépendendantes du processus $\cQ$ apparaissant dans le théorème précédent.
\end{itemize}
\end{theo}

\subsubsection{Preuve du théorème \ref{t:ABBS-normalise} à partir du théorème \ref{t:ABBS-2.3-normalise}}

Pour tout $k \ge 1$, on note $\cH_k$ l'ensemble des particules qui sont les premières de leur lignée à atteindre le niveau $k$.
Cet ensemble est presque sûrement fini.
Notons $H_k$ son cardinal et posons :
$$
Z_k = ke^{-k} H_k.
$$
Le processus $(Z_k)_k$ est la martingale $(Z(t))_t$ regardée le long de lignes d'arrêts différentes.
La suite des $Z_k$ converge également presque sûrement vers $Z$.
Si $u$ est un élément de $\cH_k$ et si $t$ est suffisamment grand, on note $X_1^{(u)}(t)$ la position de la particule la plus basse 
à l'instant $t$ parmi celles qui descendent de $u$.

Posons, pour tout $t$ assez grand :
$$
\cP^*_{k,t} = \sum_{u \in \cH_k}\delta_{X_1^{(u)}(t)} - m(t)+\ln(C_wZ_k).
$$
En utilisant la propriété de branchement du mouvement brownien branchant, le fait que $m(t+c)-m(t)$ converge vers $0$ pour tout $c$
et le fait que $X_1(t)-m(t)$ converge vers $W$, on obtient la convergence de $\cP^*_{k,t}$ vers :
$$
\cP^*_{k,\infty} = \sum_{u \in \cH_k} \delta_{W(u)} + k+\ln(C_wZ_k)
$$
où, conditionnellement à $\cF_{\cH_k}$, les $W(u)$ sont des copies indépendantes de $W$.

En utilisant \eqref{e:asymptotiquewABBS} et des résultats classiques de la théorie des valeurs extrêmes, on obtient facilement la convergence en loi de
$$
\sum_{u \in \cH_k} \delta_{W(u)} + (\ln(H_k)+\ln(\ln(H_k))+\ln(C_w))
$$
vers un processus ponctuel de Poisson $\chi$ sur $\R$ de mesure d'intensité $e^xdx$.
En utilisant la relation $H_k = k^{-1}e^kZ_k$, on obtient le comportement presque sûr suivant :
$$
\ln(H_k)+\ln(\ln(H_k)) = \ln(Z_k)+k+o_k(1).
$$
En combinant les deux derniers résultats, on obtient la convergence en loi de $\cP^*_{k,\infty}$ vers le processus de Poisson ponctuel $\chi$.

Pour tout $u\in\cH_k$, tout $\zeta>0$ et pour $t$ suffisamment grand, notons $\cQ_{t,\zeta}^{(u)}$ le processus ponctuel - vu depuis $X_1^{(u)}(t)$ - 
des descendants de $u$ qui se sont séparés de $X_1^{(u)}(t)$ après l'instant $t-\zeta$.
Définissons de la même manière $\cQ_{\zeta}$ à partir de $\cQ$ 
($\cQ$ ne contient pas d'information sur la généalogie, mais on peut définir $\cQ_{\zeta}$ en reprenant la construction de $\cQ$).

Fixons $\eta>0$.
Considérons l'évènement défini par les deux conditions suivantes :
\begin{enumerate}
\item Si $X_i(t) \le 2\eta$ et $X_j(t) \in [-2\eta,2\eta]$ alors $\tau_{i,j} \not\in [\zeta,t-\zeta]$.
\item Aucune particule n'atteint le niveau $k$ avant l'instant $\zeta$.
\item On a $|\ln(C_wZ)|\le \eta$.
\end{enumerate}
Une légère variante du théorème \ref{t:genealogie} permet de s'assurer que cet évènement est de probabilité proche de $1$ pour un bon choix des paramètres.
Si cet évènement est vérifié, on a, pour $t$ suffisamment grand :
$$
{\cN_Z}_{\big|[-\eta,\eta]} = \sum_{u \in \cH_k}\left(\cQ_{t,\zeta}^{(u)}+X_{1,t}^{(u)}-m(t)+\ln(C_wZ)\right)_{\big|[-\eta,\eta]}.
$$
Conditionnellement à $\cF_{\cH_k}$, les $\cQ_{t,\zeta}^{(u)}$ et les $X_{1,t}^{(u)}-m(t)+\ln(C_wZ_k)$ sont indépendants.
En faisant dépendre convenablement $\zeta$ de $k$, en faisant tendre $t$ vers l'infini puis $k$ vers l'infini, en utilisant la convergence de $\cP^*_{k,t}$
et de $\cP^*_{k,\infty}$ et le théorème \ref{t:ABBS-2.3-normalise} puis en faisant tendre $\eta$ vers l'infini,
on obtient la convergence en loi de $\cN_Z$ vers $\cQ$.

\subsubsection{Décomposition en épine dorsale}

Considérons $F$ une fonction raisonnable de $C([0,t],\R)$ dans $\R$.
Rappelons que $X_{k,t}$ désigne la trajectoire sur $[0,t]$ de la particule en position $X_k(t)$ à l'instant $t$.
On a :
$$
E\left(\sum_{k=1}^{N(t)} F(X_{k,t})\right) = E\big(N(t)\big) E\big(F(\sigma B(s)+2s, s \in [0,t])\big) = e^t E\big(F(\sigma B(s)+2s, s \in [0,t])\big)
$$
où $B$ est un mouvement brownien standard.
Un changement de mesure {\em via} la formule de Cameron-Martin permet d'éliminer la dérive et le facteur $e^t$.
On obtient :
$$
E\left(\sum_{k=1}^{N(t)} F(X_{k,t})\right) = E\big(F(\sigma B(s), s \in [0,t])\big).
$$

On peut en réalité faire un changement de mesure pour le mouvement brownien branchant lui-même.
Le processus sous la nouvelle mesure admet une description probabiliste agréable.
Détaillons cela.
On définit une martingale $(M(t))_{t \ge 0}$ en posant :
$$
M(t) = \sum_{k=1}^{N(t)} e^{-X_k(t)}.
$$
Cette martingale converge presque sûrement vers $0$.
Notons $(\cF_t)_t$ la filtration naturelle du mouvement brownien branchant.
Soit $Q$ la mesure de probabilité sur $\cF_\infty$ telle que, pour tout $t \ge 0$ :
$$
Q_{|\cF_t} = M(t).P_{|\cF_t}.
$$
Chauvin et Rouault ont fourni dans \cite{Chauvin-Rouault-88} la description suivante du processus sous la mesure de probabilité $Q$.
Une particule initialement en l'origine évolue suivant un mouvement brownien sans dérive et de variance $\sigma^2=2$.
Cette particule est l'épine dorsale du processus.
Avec un taux $2$, l'épine dorsale donne naissance à une particule qui évolue suivant un mouvement brownien branchant de loi $P$
(branchements à taux $1$ et mouvement brownien de dérive $\rho=2$ et de variance $\sigma^2=2$).
Par ailleurs, si $\Xi(t) \in \{1,\dots,N(t)\}$ désigne l'indice de l'épine dorsale à l'instant $t$, on a, pour $i \le N(t)$ :
$$
Q(\Xi(t)=i | \cF_t) = \frac{e^{-X_i(t)}}{M(t)}.
$$
Si par exemple $F(X_{1,t})$ est une fonction de la trajectoire $X_{1,t}(\cdot)$ sur $[0,t]$ de la particule qui réalise le minimum à l'instant $t$, 
on vérifie facilement :
\begin{eqnarray*}
E_P(F(X_{1,t})) 
 & = & E_P\left( \sum_{k=1}^{N(t)} F(X_{k,t})1_{k=1}\right), \\
 & = & E_Q\left( \frac{1}{M(t)}\sum_{k=1}^{N(t)} F(X_{k,t})1_{k=1}\right), \\
 & = & E_Q\left( e^{X_{\Xi(t)}(t)}F(X_{\Xi(t),t})1_{\Xi(t)=1}\right).
\end{eqnarray*}
Ce type de formule est connu sous le nom de {\em many-to-one formula}.
L'égalité $\Xi(t)=1$ est vérifiée lorsqu'aucun des minima des mouvements browniens branchants créés le long de la trajectoire de l'épine dorsale n'est, 
à l'instant $t$, sous $X_{\Xi(t)}(t)$.
Mais les instants de créations de ces mouvements browniens branchants forment un processus de Poisson ponctuel sur $[0,t]$ de mesure d'intensité $2dt$.
Par ailleurs, si $s$ est l'un de ces instants, la probabilité que son minimum soit à l'instant $t$ sous $X_{\Xi(t)}(t)$ est (conditionnellement à la trajectoire 
de l'épine dorsale et aux temps de branchement sur cette épine dorsale) :
$$
G_{t-s}(X_{\Xi(t)}(t)-X_{\Xi(t)}(s))
$$
où, pour tout $t \ge 0$ et tout $x \in \R$ :
$$
G_t(x) = P(X_1(t) \le x).
$$
Ainsi, le nombre de mouvements browniens branchants créés dont le minimum est sous $X_{\Xi(t)}(t)$ à l'instant $t$ suit, conditionnellement à la trajectoire de l'épine dorsale,
une loi de Poisson de paramètre $2\int_0^t G_{t-s}(X_{\Xi(t)}(t)-X_{\Xi(t)}(s)) ds$.
Par conséquent :
$$
E_P(F(X_{1,t})) = E_Q\left( e^{X_{\Xi(t)}(t)}F(X_{\Xi(t),t})e^{-2\int_0^t G_{t-s}(X_{\Xi(t)}(t)-X_{\Xi(t)}(s)) ds}\right).
$$
Notons $(B(t))_t$ un mouvement brownien standard. La trajectoire de l'épine dorsale étant un mouvement brownien sans dérive de variance $\sigma^2$ on obtient finalement :
\begin{equation} \label{e:manytoone}
E_P(F(X_{1,t})) = E\left( e^{\sigma B(t)}F(\sigma B(s), s \in [0,t])e^{-2\int_0^t G_{t-s}(\sigma B(t) - \sigma B(s)) ds}\right).
\end{equation}

\subsubsection{Preuve du théorème \ref{t:ABBS-2.3-normalise}}

Soit $F_1$ une fonction raisonnable de $C(\R_+,\R)$ dans $\R$.
Si $X$ est une fonction continue de $[0,t]$ dans $\R$, $F_1(X)$ est l'image par $F_1$ de la fonction définie par $s \mapsto X(\min(t,s))$. 
Soient $A_1, \dots, A_n$ des sous-ensembles boréliens de $\R$ et $\alpha_1, \dots, \alpha_n$ des nombres réels.
Soit enfin $\eta > 0$.
Nous nous intéressons au comportement asymptotique de $I(x,t,\zeta)$ défini pour un paramètre $x >0$ et pour $t \ge \zeta \ge 0$ par :
$$
I(x,t,\zeta) = E\left(1_A(X_{1,t}) F_1(Y_t(s), s \in [0,\zeta])) e^{-\sum_{i=1}^n \alpha_i \cQ(t,\zeta)(A_i)}1_{|X_1(t)-m(t)| \le \eta}\right).
$$
Une trajectoire $X:[0,t]\to\R$ appartient à $A$ si : 
\begin{itemize}
\item On a $X(s) \ge a(s)$ pour $s \in [0,t-\zeta]$ où $a(s)=-x$ pour $s \le t/2$ et $a(s)=m(t)-x$ pour $s > t/2$.
\item On a $X(s)-X(t) \ge -x$ pour $s \in [t-\zeta,t]$.
\item On a $X(t-\zeta)-X(t) \in [\zeta^{1/3},\zeta^{2/3}]$.
\item Le minimum de $X$ sur $[t/2,t]$ est atteint en un point de $[t-x,t]$.
\end{itemize}
Avec une probabilité abitrairement proche de $1$ pour $x, \zeta$ et $t$ suffisamment grands, on a 
$X_1(t)-m(t) \in [-\eta,\eta]$ ou $X_{1,t} \in  A$. 
C'est un résultat semblable au théorème \ref{t:localisation}.

Notons $E_{a,b}^{(t)}$ l'espérance sous laquelle $B=(B(s))_{s \in [0,t]}$ est un pont brownien de longueur $t$ entre $a$ et $b$.
Posons 
$$
G_s^*(a) = 1 - E\left(e^{-\sum_{i=1}^n \alpha_i 1_{A_i}(\cN(s)-a)}1_{\min \cN(s) -a \ge 0}\right).
$$
Notons que, si tous les $\alpha_i$ sont nuls, on a $G_s^*(a) = P(\min \cN(s)-a < 0)=G_s(a)$.
Un raisonnement semblable à celui menant à \eqref{e:manytoone},
un renversement de trajectoire - remplacement de $B=(B(s))_{s \in [0,t]}$ par le processus de même loi $\overline{B}=(B(t)-B(t-s))_{s \in [0,t]}$) - 
puis un conditionnement/déconditionnement par la valeur de $\sigma B(t)$ permet d'obtenir :
$$
I(x,t,\zeta) = \frac{C}{t} \int_{-\eta}^{\eta} P(\sigma B(t) \in m(t)+dz)e^{m(t)+z} J(x,t,\zeta,z)
$$
où :
$$
J(x,t,\zeta,z) = t E_{0,\sigma^{-1}(m(t)+z}^{(t)}\left(1_A( \sigma \overline{B}) F_1(\sigma B)
e^{-2\int_0^{\zeta} G^*_s(\sigma B_s)ds-2\int_{\zeta}^t G_s(\sigma B_s)ds}\right).
$$
En explicitant la densité de la loi de $\sigma B(t)$ et en exploitant le fait que les $z$ pertinents appartiennent à un intervalle borné, on obtient :
\begin{equation}\label{e:IJ}
I(x,t,\zeta) \sim_{t\to\infty} \frac{C}{t^{3/2}} \int_{-\eta}^{\eta} e^{m(t)+z} J(x,t,\zeta,z) dz = C \int_{-\eta}^{\eta} e^z J(x,t,\zeta,z) dz
\end{equation}
où $C$ est une constante qui change de ligne en ligne.
Il reste à étudier le comportement asymptotique de $J(x,t,\zeta,z)$ pour $z$ dans $[-\eta,\eta]$.
Pour cela, on applique la propriété de Markov au pont brownien à l'instant $\zeta$. 
Rappelons que nous avons renversé la trajectoire, cela correspond donc au temps réel $t-\zeta$.
Notons $\theta$ l'instant où $(B(s))_{s \in [0,\zeta]}$ atteint son maximum.
On obtient, toujours en exploitant le fait que $z$ appartient à un intervalle borné :
$$
J(x,t,\zeta,z) \sim_{t \to \infty} \int_{-\zeta^{2/3}}^{-\zeta^{1/3}} K(x,\zeta,w)L(x,t,\zeta, z,w)P(\sigma B_{\zeta} \in dw),
$$
où
$$
K(x,\zeta,w) = E_{0,\sigma^{-1}w}^{(\zeta)}\left(1_{\sigma B(s) \le x, s \in [0,\zeta]}1_{\theta \le x} 
F_1(\sigma B(s), s \in [0,\zeta])e^{-2\int_0^{\zeta} G^*_s(\sigma B_s)ds}\right)
$$
et où
$$
L(x,t,\zeta, z,w) = t E_{0,\sigma^{-1}(m(t)+z-w)}^{(t-\zeta)} \left(1_{\sigma B(t-\zeta)-\sigma B(t-\zeta-s) \ge a(s), s \in [0,t-\zeta]}
e^{-2\int_{0}^{t-\zeta} G_{s+\zeta}(w+\sigma B_s)ds}\right).
$$

\ABBetS montrent alors que $L(x,t,\zeta,z,w)$ converge lorsque $t$ tend vers l'infini et que cette limite vérifie, uniformément en $\zeta$ :
\begin{equation}\label{e:6.1}
\lim_{t \to \infty} L(x,t,\zeta,z,w) \sim_{w \to -\infty} |w| \varphi(x,z)
\end{equation}
pour une fonction $\varphi$ explicite.
On peut alors en déduire, en utilisant le fait que les $w$ pertinents appartiennent à l'intervalle $[-\zeta^{2/3},-\zeta^{1/3}]$ :
$$
\lim_{t\to\infty} J(x,t,\zeta,z)= \int_{-\zeta^{2/3}}^{-\zeta^{1/3}} K(x,\zeta,w)|w|\varphi(x,z)P(\sigma B_{\zeta} \in dw).
$$
On peut alors montrer :
\begin{equation}\label{e:6.2}
\lim_{\zeta \to \infty} \lim_{t \to \infty} J(x,t,\zeta,z) = \varphi(x,z) \sigma \int_0^{\sigma^{-1}x} 
E\left(F_1(\sigma \Gamma^{(b)}(s), s \ge 0) e^{-2\int_0^{\infty} G^*_s(\sigma \Gamma_s^{(b)})ds} 1_{T_b \le x} \right)db.
\end{equation}
On conlut la preuve en revenant à \eqref{e:IJ} et laissant tendre $x$ vers l'infini.

La relation \eqref{e:6.2} s'obtient par des décompositions de trajectoires du mouvement brownien.
L'obtention de la relation \eqref{e:6.1} est plus délicate.
Des arguments semblables à ceux menant au théorème \ref{t:genealogie} sur la généalogie permettent de limiter l'intégrale sur $[0,t-\zeta]$
intervant dans l'exponentielle à une intégrale sur des ensembles de la forme $[0,M] \cup [t-\zeta-M,t-\zeta]$.
Lorsque $w$ tend vers $-\infty$, l'intégrale sur $[0,M]$ devient triviale.
Après retournement du temps, l'intégrale sur $[t-\zeta-M,t-\zeta]$ s'écrit :
$$
\int_0^M G_{t-s}(m(t)+z+\sigma B^{(t)}(s) )ds
$$
où la loi de $B^{(t)}(s)$ est, asymptotiquement en $t$, celle d'une mouvement brownien standard. La convergence de l'intégrale se déduit alors
de la convergence en loi de $X_1(t)-m(t)$.
Nous renvoyons aux sections 7, 8 et 9 de \cite{ABBS} pour les détails des preuves de \eqref{e:6.1} et \eqref{e:6.2}.

\bibliographystyle{alpha}
\bibliography{bibBourbaki}

\end{document}